\documentclass[a4paper,12pt]{amsart}

\usepackage{amssymb,amsthm,amsmath,amscd,amsfonts,bbm,mathrsfs}

\usepackage[cmtip,all]{xy}

\theoremstyle{plain}
\newtheorem{Lemma}{Lemma}
\newtheorem{Thm}[Lemma]{Theorem}
\newtheorem*{Thm*}{Theorem}
\newtheorem{Prop}[Lemma]{Proposition}
\newtheorem{Cor}[Lemma]{Corollary}
\theoremstyle{definition}
\newtheorem{Defn}[Lemma]{Definition}

\newtheorem*{Example*}{Example}
\theoremstyle{remark}
\newtheorem{Remark}[Lemma]{Remark}
\newtheorem*{Remark*}{Remark}
\newtheorem{Remarks}[Lemma]{Remarks}

\numberwithin{Lemma}{section}
\numberwithin{equation}{section}

\newcommand{\BBB}{\mathscr{B}}

\newcommand{\DDD}{\mathscr{D}}
\newcommand{\FFF}{\mathscr{F}}
\newcommand{\HHH}{\mathscr{H}}

\newcommand{\NNN}{\mathscr{N}}
\newcommand{\PPP}{\mathscr{P}}
\newcommand{\WWW}{\mathscr{W}}

\newcommand{\Fa}{\mathfrak{a}}
\newcommand{\Fb}{\mathfrak{b}}
\newcommand{\Fc}{\mathfrak{c}}
\newcommand{\Fm}{\mathfrak{m}}
\newcommand{\Fn}{\mathfrak{n}}
\newcommand{\Fp}{\mathfrak{p}}
\newcommand{\FS}{\mathfrak{S}}

\newcommand{\II}{{\mathbb{I}}}
\newcommand{\MM}{{\mathbb{M}}}
\newcommand{\NN}{{\mathbb{N}}}
\newcommand{\QQ}{{\mathbb{Q}}}
\newcommand{\WW}{{\mathbb{W}}}
\newcommand{\ZZ}{{\mathbb{Z}}}

\DeclareMathOperator{\id}{id}

\DeclareMathOperator{\Hom}{Hom}
\DeclareMathOperator{\Ext}{Ext}
\DeclareMathOperator{\Ker}{Ker}
\DeclareMathOperator{\Coker}{Coker}
\DeclareMathOperator{\Spec}{Spec}

\DeclareMathOperator{\Rad}{Rad}

\newcommand{\tvarkappa}{\tilde\varkappa}
\newcommand{\tvarphi}{\tilde\varphi}
\newcommand{\tpsi}{\tilde\psi}
\newcommand{\tsigma}{\tilde\sigma}
\newcommand{\tBBB}{\tilde\BBB}

\newcommand{\tI}{\tilde I}
\newcommand{\tf}{\tilde f}

\begin{document}

\title[Frames and finite group schemes]
{Frames and finite group schemes over
complete regular local rings}
\author{Eike Lau}
\date{\today}
\address{Fakult\"{a}t f\"{u}r Mathematik,
Universit\"{a}t Bielefeld, D-33501 Bielefeld, Germany}
\email{lau@math.uni-bielefeld.de}
% \subjclass[2000]{Primary 14L05; Secondary 14F30}

\begin{abstract}
Let $p$ be an odd prime.
We show that the classification of $p$-divisible groups
by Breuil windows and the classification of finite flat
group schemes of $p$-power order by Breuil modules
hold over any complete regular local ring with perfect
residue field of characteristic $p$.
We use a formalism of frames and windows
with an abstract deformation theory
that applies to Breuil windows.
\end{abstract}

\maketitle

\section{Introduction}

Let $R$ be a complete regular local ring with 
perfect residue field $k$ of odd characteristic $p$.
One can write $R=\FS/E\FS$ with
$$
\FS=W(k)[[x_1,\ldots,x_r]]
$$
such that $E\in\FS$ is a power series with constant term $p$. 
Let $\sigma$ be the continuous endomorphism of $\FS$
that extends the Frobenius automorphism of $W(k)$ by
$\sigma(x_i)=x_i^p$. 
Following Vasiu and Zink, a {\em Breuil window} 
relative to $\FS\to R$ is a pair $(Q,\phi)$ 
where $Q$ is a free $\FS$-module of finite rank, 
and where $$\phi:Q\to Q^{(\sigma)}$$ is an $\FS$-linear 
homomorphism with cokernel annihilated by $E$.

\begin{Thm}
\label{Th-equiv}
The category of $p$-divisible groups over $R$ is equivalent
to the category of Breuil windows relative to $\FS\to R$.
\end{Thm}
 
If $R$ has characteristic $p$ this follows from more general
results of A.~de Jong \cite{deJong}; this case is included here
only for completeness.
If $r=1$ and $E$ is an Eisenstein polynomial, 
Theorem \ref{Th-equiv} was conjectured by Breuil \cite{Breuil} 
and proved by Kisin \cite{Kisin}. 
When $E$ is the deformation of an Eisenstein polynomial
the result is proved in \cite{Vasiu-Zink-Classif}. 

Like in these cases one can derive 
a classification of finite group schemes:
A {\em Breuil module} relative to $\FS\to R$ is a triple
$(M,\varphi,\psi)$ where $M$ is a finitely generated 
$\FS$-module annihilated by a power of $p$ and of 
projective dimension at most one, and where
$$
\varphi:M\to M^{(\sigma)},\qquad\psi:M^{(\sigma)}\to M
$$
are homomorphisms of $\FS$-modules with $\varphi\psi=E$
and $\varphi\psi=E$. If $R$ has characteristic zero
such triples are equivalent to pairs $(M,\varphi)$ 
such that the cokernel of $\varphi$ is annihilated by $E$.
 
\begin{Thm}
\label{Th-equiv-fin}
The category of finite
flat group schemes over $R$ annihilated by a power of $p$
is equivalent to the category of Breuil modules relative
to $\FS\to R$. 
\end{Thm} 

This result is applied in \cite{Vasiu-Zink-Finite} to
the question whether abelian schemes or $p$-divisible
groups defined over $\Spec R\setminus\{\Fm_R\}$ extend to
$\Spec R$.

\subsection*{Frames and windows}

To prove Theorem \ref{Th-equiv} we show that Breuil windows 
are equivalent to Dieudonn\'e displays over $R$, which are 
equivalent to $p$-divisible groups by \cite{Zink-DDisp}; 
the same route is followed in \cite{Vasiu-Zink-Classif}.
So the main part of this article is purely module theoretic:

We introduce a notion of frames and windows 
(motivated by \cite{Zink-Windows}) which allows
to formulate a deformation theory that generalises 
the deformation theory of Dieudonn\'e displays 
and that also applies to Breuil windows. 
Technically the main point is the formalism of $\sigma_1$;
the central result is the lifting of windows
in Theorem \ref{Th-crys}.

This is applied as follows.
For each positive integer $a$ we consider the rings 
$\FS_a=\FS/(x_1,\ldots,x_r)^a\FS$ and $R_a=R/\Fm_R^a$.
There is an obvious notion of Breuil windows
relative to $\FS_a\to R_a$ and a functor
$$
\varkappa_a:(\text{Breuil windows rel.\ $\FS_a\to R_a$})
\to(\text{Dieudonn\'e displays}/R_a).
$$
The deformation theory implies that on both sides 
lifts from $a$ to $a+1$ are classified by lifts 
of the Hodge filtration in a compatible way. 
Thus $\varkappa_a$ is an equivalence for all $a$
by induction, and Theorem \ref{Th-equiv} follows.

\subsection*{Complements}

There is some freedom in the choice of 
the Frobenius lift on $\FS$. Namely,
let $\sigma$ be a ring endomorphism of $\FS$
which preserves the ideal $J=(x_1,\ldots,x_r)$ 
and which induces the Frobenius on $\FS/p\FS$.
If the endomorphism $\sigma/p$ of $J/J^2$ is
nilpotent modulo $p$, Theorems \ref{Th-equiv}
and \ref{Th-equiv-fin} hold without change.

All of the above equivalences of categories are compatible 
with the natural duality operations on both sides.

If the residue field $k$ is not assumed perfect there is 
an analogue of Theorems \ref{Th-equiv} and \ref{Th-equiv-fin}
for connected groups. Here $p=2$ is allowed. 
The ring $W(k)$ is replaced by a $p$-ring of $k$, 
and the operators $\phi$ and $\varphi$ must be 
nilpotent modulo the maximal ideal of $\FS$. 

In the first version of this article \cite{Lau-Win}
the formalism of frames was introduced only to give
an alternative proof of the results of 
Vasiu and Zink \cite{Vasiu-Zink-Classif}. In response, 
they pointed out that both their and this approach 
apply in more generality, e.g.\ in the case where 
$E\in\FS$ takes the form $E=g+p\epsilon$ such
that $\epsilon$ is a unit and $g$ divides $\sigma(g)$.
However, the method of loc.~cit.\ seems not to
give Theorem \ref{Th-equiv} completely.

\medskip

All rings in this text are commutative and have a unit.

\medskip
\noindent
{\it Acknowledgements.}
The author thanks A.~Vasiu and Th.~Zink for valuable 
discussions, and in particular Th.~Zink for sharing 
his notion of $\varkappa$-frames
and for suggesting to include sections
\ref{Se-disp} and \ref{Se-gen}.

%------------------------------------------------------------
%\cleardoublepage
\section{Frames and windows}
\label{Se-fr-win}

Let $p$ be a prime.
The following notion of frames and windows 
differs from \cite{Zink-Windows}. 
Some definitions and arguments could be
simplified by assuming that the relevant
rings are local, which is the case in our
applications, but we work in more generality 
until section \ref{Se-abs-def}.

If $\sigma:S\to S$ is a ring endomorphism, 
for an $S$-module $M$ we write 
$M^{(\sigma)}=S\otimes_{\sigma,S}M$, and for a $\sigma$-linear
homomorphism $g:M\to N$ we denote
by $g^{\sharp}:M^{(\sigma)}\to N$ its linearisation,
$g^{\sharp}(s\otimes m)=sg(m)$.

\begin{Defn}
\label{Def-frame}
A frame is a quintuple 
$$
\FFF=(S,I,R,\sigma,\sigma_1)
$$
consisting of a ring $S$, an ideal $I$ of $S$,
the quotient ring $R=S/I$,
a ring endomorphism $\sigma:S\to S$, and a $\sigma$-linear
homomorphism of $S$-modules $\sigma_1:I\to S$, such that
the following conditions hold.
\begin{enumerate}
\renewcommand{\theenumi}{\roman{enumi}}
\item\label{Def-frame-Rad}
$I+pS\subseteq\Rad(S)$,
\item\label{Def-frame-Frob}
$\sigma(a)\equiv a^p$ mod $pS$ for $a\in S$,
\item\label{Def-frame-surj}
$\sigma_1(I)$ generates $S$ as an $S$-module.
\newcounter{frame-count}
\setcounter{frame-count}{\value{enumi}}
\end{enumerate}
\end{Defn}

\begin{Remark*}
With some modifications the theory also works without 
assuming \eqref{Def-frame-surj}; see section \ref{Se-gen}.
In our examples $\sigma_1(I)$ contains $1$.
\end{Remark*}

\begin{Lemma}
\label{Le-theta}
For every frame $\FFF$ there is a unique element
$\theta\in S$ such that $\sigma(a)=\theta\sigma_1(a)$
for $a\in I$. 
\end{Lemma}

\begin{proof}
Condition \eqref{Def-frame-surj} means that 
$\sigma_1^\sharp:I^{(\sigma)}\to S$ is surjective. 
If $b\in I^{(\sigma)}$ satisfies
$\sigma_1^\sharp(b)=1$, then necessarily 
$\theta=\sigma^\sharp(b)$. For $a\in I$ we compute 
$\sigma(a)=\sigma_1^\sharp(b)\sigma(a)=
\sigma_1^\sharp(ba)=\sigma^\sharp(b)\sigma_1(a)$
as desired.
\end{proof}

\begin{Defn}
\label{Def-win}
A window over a frame $\FFF$ is a quadruple 
$$
\PPP=(P,Q,F,F_1)
$$ 
where $P$ is a finitely generated projective
$S$-module, $Q\subseteq P$ is a submodule, 
$F:P\to P$ and $F_1:Q\to P$ are $\sigma$-linear
homomorphisms of $S$-modules, such that
the following conditions hold.
\begin{enumerate}
\item \label{Def-win-norm}
There is a decomposition $P=L\oplus T$ with $Q=L\oplus IT$,
\item \label{Def-win-FF1}
$F_1(ax)=\sigma_1(a)F(x)$ for $a\in I$ and $x\in P$, 
\item \label{Def-win-surj}
$F_1(Q)$ generates $P$ as an $S$-module.
\end{enumerate}
A decomposition as in \eqref{Def-win-norm} 
is called a normal decomposition.
\end{Defn}

\begin{Remark*}
The operator $F$ is determined by $F_1$. 
Indeed, if $b\in I^{(\sigma)}$ satisfies
$\sigma_1^\sharp(b)=1$, then condition \eqref{Def-win-FF1} 
implies that $F(x)=F_1^\sharp(bx)$ for $x\in P$. 
In particular we have $F(x)=\theta F_1(x)$ when $x$ lies in $Q$.
\end{Remark*}

\begin{Remark}
\label{Re-lift}
Condition \eqref{Def-win-norm} implies that
\begin{enumerate}
\item[(\ref{Def-win-norm}')]
$P/Q$ is a projective $R$-module.
\end{enumerate}
If finitely generated projective $R$-modules lift 
to projective $S$-modules, necessarily finitely 
generated because $I\subseteq\Rad(S)$,
then \eqref{Def-win-norm} is equivalent to
(\ref{Def-win-norm}'). In all our examples,
this lifting property holds because 
$S$ is local or $I$-adically complete.
\end{Remark}

We recall that a $\sigma$-linear isomorphism is a 
$\sigma$-linear homomorphism with bijective linearisation.

\begin{Lemma}
\label{Le-norm}
Let $\FFF$ be a frame,
let $P=L\oplus T$ be a finitely generated projective
$S$-module, and let $Q=L\oplus IT$.
The set of $\FFF$-window structures $(P,Q,F,F_1)$ 
on these modules is mapped bijectively to the set 
of $\sigma$-linear isomorphisms 
$$
\Psi:L\oplus T\to P
$$
by the assignment $\Psi(l+t)=F_1(l)+F(t)$
for $l\in L$ and $t\in T$.
\end{Lemma}

The triple $(L,T,\Psi)$ is called a normal representation
of $(P,Q,F,F_1)$.

\begin{proof}
If $(P,Q,F,F_1)$ is an $\FFF$-window, by \eqref{Def-win-FF1}
and \eqref{Def-win-surj} the linearisation
of the associated homomorphism $\Psi$ is surjective, 
thus bijective as $P$ and $P^{(\sigma)}$ 
are projective $S$-modules of equal rank by  
\eqref{Def-frame-Rad} and \eqref{Def-frame-Frob}. 
Conversely, if $\Psi$ is given, one gets
an $\FFF$-window by $F(l+t)=\theta\Psi(l)+\Psi(t)$
and $F_1(l+at)=\Psi(l)+f_1(a)\Psi(t)$ for 
$l\in L$, $t\in T$, and $a\in I$.
\end{proof}

\begin{Example*}
The Witt frame of a $p$-adically complete ring $R$ is
$$
\WWW_R=(W(R),I_R,R,f,f_1)
$$
where $f$ is the Frobenius endomorphism and where
$f_1:I_R\to W(R)$ is the inverse of the Verschiebung
homomorphism. Here $\theta=p$.
We have $I_R\subset\Rad(W(R))$ because $W(R)$ 
is $I_R$-adically complete; see \cite[Proposition 3]{Zink-Disp}.
Windows over $\WWW_R$ are 3n-displays over $R$ in the 
sense of \cite{Zink-Disp}, called displays 
in \cite{Messing-Bour}, which is the terminology we follow.
\end{Example*}

%------------------------------------------------------------
\subsection*{Functoriality}

Let $\FFF$ and $\FFF'$ be frames.

\begin{Defn}
A homomorphism of frames $\alpha:\FFF\to\FFF'$ is 
a ring homomorphism $\alpha:S\to S'$ with
$\alpha(I)\subseteq I'$ such that 
$\sigma'\alpha=\alpha\sigma$ and 
$\sigma_1'\alpha=u\cdot\alpha\sigma_1$
for a unit $u\in S'$. If $u=1$ then $\alpha$ is called strict.
\end{Defn}

\begin{Remark}
\label{Re-u-unique}
The unit $u$ is unique because $\alpha\sigma_1(I)$ 
generates $S'$ as an $S'$-module. If we want to
specify $u$ we say that $\alpha$ is a $u$-homomorphism. 
We have $\alpha(\theta)=u\theta'$. There is a unique
factorisation $\alpha=\omega\alpha'$ such that
$\alpha':\FFF\to\FFF''$ is strict and $\omega:\FFF''\to\FFF'$
is invertible.
\end{Remark}

Let $\alpha:\FFF\to\FFF'$ be a $u$-homomorphism of frames.

\begin{Defn}
Let $\PPP$ and $\PPP'$ be windows over $\FFF$ and
$\FFF'$, respectively.
An $\alpha$-homomorphism of windows $g:\PPP\to\PPP'$ 
is a homomorphism of $S$-modules $g:P\to P'$ with 
$g(Q)\subseteq Q'$ such that $F'g=gF$ and $F_1'g=u\cdot gF_1$.
A homomorphism of windows over $\FFF$ is an 
$\id_{\PPP}$-homomorphism in the previous sense.
\end{Defn}

\begin{Lemma}
For each window $\PPP$ over $\FFF$ there is a base change
window $\alpha_*\PPP$ over $\FFF'$ together
with an $\alpha$-homomorphism $\PPP\to\alpha_*\PPP$ 
that induces a bijection
$\Hom_{\FFF'}(\alpha_*\PPP,\PPP')=\Hom_\alpha(\PPP,\PPP')$
for all windows $\PPP'$ over $\FFF'$.
\end{Lemma}

\begin{proof}
Clearly this requirement determines $\alpha_*\PPP$ 
uniquely. It can be constructed explicitly as follows: 
If $(L,T,\Psi)$ is a normal representation of $\PPP$,
then a normal representation of $\alpha_*\PPP$ is
$(S'\otimes_SL,S'\otimes_ST,\Psi')$
with $\Psi'(s'\otimes l)=u\sigma'(s')\otimes\Psi(l)$
and $\Psi'(s'\otimes t)=\sigma'(s')\otimes\Psi(t)$.
\end{proof}

\begin{Remark*}
If $\alpha_*\PPP=(P',Q',F',F_1')$, then
$P'=S'\otimes_SP$, and $Q'$ is the image of
$S'\otimes_SQ\to P'$, which may differ from
$S'\otimes_SQ$.
\end{Remark*}

%----------------------------------------------------------------
\subsection*{Limits}

Windows are compatible with projective limits of
frames in the following sense. Assume that for each
positive integer $n$ we have a frame 
$\FFF_n=(S_n,I_n,R_n,\sigma_n,\sigma_{1n})$ and a strict homomorphism
of frames $\pi_n:\FFF_{n+1}\to\FFF_n$ such that the maps 
$S_{n+1}\to S_n$ and $I_{n+1}\to I_n$ are surjective 
and $\Ker(\pi_n)$ is contained in $\Rad(S_{n+1})$.
We obtain a frame $\varprojlim\FFF_n=(S,I,R,\sigma,\sigma_1)$
with $S=\varprojlim S_n$ etc.
By definition, a window over $\FFF_*$ is a system
$\PPP_*$ of windows $\PPP_n$ over $\FFF_n$ together
with isomorphisms $\pi_{n*}\FFF_{n+1}\cong\FFF_n$.

\begin{Lemma}
\label{Le-lim}
The category of windows over $\varprojlim\FFF_n$ is
equivalent to the category of windows over $\FFF_*$.
\end{Lemma}

\begin{proof}
The obvious functor from windows
over $\varprojlim\FFF_n$ to windows over $\FFF_*$ 
is fully faithful. We must show that for a window 
$\PPP_*$ over $\FFF_*$, the projective limit 
$\varprojlim\PPP_n=(P,Q,F,F_1)$ defined by
$P=\varprojlim P_n$ etc.\ is a window over 
$\varprojlim\FFF_n$. The condition 
$\Ker(\pi_n)\subseteq\Rad(S_{n+1})$ implies that 
$P$ is a finitely generated projective $S$-module 
and that $P/Q$ is projective over $R$.
In order that $P$ has a normal decomposition it suffices
to show that any normal decomposition of $\PPP_n$ lifts
to a normal decomposition of $\PPP_{n+1}$. Assume that
$P_n=L_n'\oplus T_n'$ and $P_{n+1}=L_{n+1}\oplus T_{n+1}$
are normal decompositions and let $P_n=L_n\oplus T_n$
be induced by the second. Since 
$T_n\otimes R_n\cong P_n/Q_n\cong T'_n\otimes R_n$ and
$L_n\otimes R_n\cong Q_n/IP_n\cong L_n'\otimes R_n$
we have $T_n\cong T'_n$ and $L_n\cong L'_n$.
Hence the two decompositions of $P_n$ differ by
an automorphism of $L_n\oplus T_n$
of the type $u=\left(\begin{smallmatrix} a & b \\ c & d
\end{smallmatrix}\right)$ with $c:L_n\to I_nT_n$. Now $u$
lifts to an endomorphism $u'=\left(\begin{smallmatrix} 
a' & b' \\ c' & d' \end{smallmatrix}\right)$ of
$L_{n+1}\oplus T_{n+1}$ with $c':L_{n+1}\to I_{n+1}T_{n+1}$, 
and $u'$ is an automorphism as 
$\Ker(\pi_n)\subseteq\Rad(S_{n+1})$.
The required lifting of normal decompositions follows.
All remaining window axioms for $\varprojlim\PPP_n$ 
are easily checked.
\end{proof}

\begin{Remark}
\label{Re-lim}
Assume that $S_1$ is a local ring. Then all $S_n$ and
$S$ are local too. Hence $\varprojlim\FFF_n$
satisfies the lifting property of Remark \ref{Re-lift},
so the normal decomposition of $P$ in the preceding
proof is automatic.
\end{Remark}

\subsection*{Duality}
Let $\PPP$ be a window over a frame $\FFF$.
The dual window $\PPP^t=(P',Q',F',F_1')$ is defined as follows.
We have $P'=\Hom_S(P,S)$ and $Q'=\{x'\in P'\mid x'(Q)\subseteq I\}$.
The operator $F_1':Q'\to P'$ is defined by the relation
$$
F_1'(x')(F_1(x))=\sigma_1(x'(x))
$$
for $x'\in Q'$ and $x\in Q$. This determines
$F_1'$ and $F'$ uniquely. If $(L,T,\Psi)$ is a normal 
representation for $\PPP$, then a normal representation 
for $\PPP^t$ is given by $(T^\vee,L^\vee,\Psi')$ where 
$(\Psi')^\sharp$ is equal to $((\Psi^\sharp)^{-1})^\vee$.
This shows that $F_1'$ and $F'$ are well-defined. 
For a more detailed exposition of the duality formalism 
in the case of (Diedonn\'e) displays we refer to 
\cite[Definition 19]{Zink-Disp} or  \cite[section 3]{Lau-Dual}. 
There is a natural isomorphism $\PPP^{tt}\cong\PPP$.
For a homomorphism of frames $\alpha:\FFF\to\FFF'$ we have a 
natural isomorphism $(\alpha_*\PPP)^t\cong\alpha_*(\PPP^t)$.

%----------------------------------------------------------------
\section{Crystalline homomorphisms}

\begin{Defn}
\label{Def-crys}
A homomorphism of frames $\alpha:\FFF\to\FFF'$ is
called crystalline if the functor
$
\alpha_*:(\text{windows over $\FFF$})
\to(\text{windows over $\FFF'$})
$
is an equivalence of categories.
\end{Defn}

\begin{Thm}
\label{Th-crys}
Let $\alpha:\FFF\to\FFF'$ be a strict homomorphism of frames that
induces an isomorphism $R\cong R'$ and a surjection $S\to S'$ 
with kernel $\Fa\subset S$. We assume that there is a finite
filtration $\Fa=\Fa_0\supseteq\ldots\supseteq\Fa_n=0$
with $\sigma(\Fa_i)\subseteq\Fa_{i+1}$ and 
$\sigma_1(\Fa_i)\subseteq\Fa_i$ such that $\sigma_1$ is 
elementwise nilpotent on $\Fa_i/\Fa_{i+1}$.
We assume that finitely generated projective $S'$-modules
lift to projective $S$-modules. Then $\alpha$ is crystalline.
\end{Thm}

In many applications the lifting property of projective
modules holds because $\Fa$ is nilpotent or $S$ is local. 
The proof of Theorem \ref{Th-crys}
is a variation of the proofs of \cite[Theorem 44]{Zink-Disp} 
and \cite[Theorem 3]{Zink-DDisp}. 

\begin{proof}
The homomorphism $\alpha$ factors into 
$\FFF\to\FFF''\to\FFF'$ where the frame $\FFF''$ 
is determined by $S''=S/\Fa_1$, so
by induction we may assume that $\sigma(\Fa)=0$. 
The functor $\alpha_*$ is essentially surjective 
because normal representations $(L,T,\Psi)$
can be lifted from $\FFF'$ to $\FFF$.
In order that $\alpha_*$ is fully faithful
it suffices that $\alpha_*$ is fully faithful on 
automorphisms because a homomorphism $g:\PPP\to\PPP'$ 
can be encoded by the automorphism 
$\left(\begin{smallmatrix}1&0\\g&1\end{smallmatrix}\right)$
of $\PPP\oplus\PPP'$.
Since for a window $\PPP$ over $\FFF$ 
an automorphism of $\alpha_*\PPP$
can be lifted to an $S$-module automorphism of $P$ 
it suffices to prove the following assertion.

{\em
Assume that $\PPP=(P,Q,F,F_1)$ and $\PPP'=(P,Q,F',F_1')$
are two windows over $\FFF$ such that $F\equiv F'$ and
$F_1\equiv F_1'$ modulo $\Fa$. Then there is a unique
isomorphism $g:\PPP\cong\PPP'$ with 
$g\equiv\id$ modulo $\Fa$.
}

We write $F_1'=F_1+\eta$ and $F'=F+\varepsilon$ and $g=1+\omega$, 
where the $\sigma$-linear homomorphisms $\eta:Q\to\Fa P$ and 
$\varepsilon:P\to\Fa P$ are given, and where $\omega:P\to\Fa P$ 
is an arbitrary homomorphism of $S$-modules.
The induced $g$ is an isomorphism of windows if and only if
$gF_1=F_1'g$ on $Q$, which translates into
\begin{equation}
\label{EqCond1}
\eta=\omega F_1-F'_1\omega.
\end{equation}
We fix a normal decomposition $P=L\oplus T$, thus $Q=L\oplus IT$.
For $l\in L$, $t\in T$, and $a\in I$ we have
\begin{gather*}
\eta(l+at)=\eta(l)+\sigma_1(a)\varepsilon(t), \\
\omega(F_1(l+at))=\omega(F_1(l))+\sigma_1(a)\omega(F(t)), \\
F_1'(\omega(l+at))=F_1'(\omega(l))+\sigma_1(a)F'(\omega(t)).
\end{gather*}
Here $F'\omega=0$ because for $a\in\Fa$ and $x\in P$
we have $F'(ax)=\sigma(a)F'(x)$, and $\sigma(\Fa)=0$.
As $\sigma_1(I)$ generates $S$ we see that
\eqref{EqCond1} is equivalent to:
\begin{equation}
\label{EqCond2}
\left\{
\begin{array}{l}
\varepsilon=\omega F\qquad\qquad\;\text{ on } T, \\
\eta=\omega F_1-F_1'\omega\quad\text{ on } L.
\end{array}
\right.
\end{equation}

As $\Psi:L\oplus T\xrightarrow{F_1+F} P$ is a $\sigma$-linear
isomorphism, the datum of $\omega$ is equivalent to 
the pair of $\sigma$-linear homomorphisms
$$
\omega_L=\omega F_1:L\to\Fa P, \qquad
\omega_T=\omega F:T\to\Fa P.
$$
Let $\lambda:L\to L^{(\sigma)}$ be the composition
$L\subseteq P\xrightarrow{(\Psi^\sharp)^{-1}}
L^{(\sigma)}\oplus T^{(\sigma)}\xrightarrow{pr_1}L^{(\sigma)}$
and let $\tau:L\to T^{(\sigma)}$ be analogous
with $pr_2$ in place of $pr_1$. 
Then the restriction $\omega|_L$ is equal to 
$\omega_L^\sharp\lambda+\omega_T^\sharp\tau$,
and \eqref{EqCond2} becomes:
\begin{equation}
\label{EqCond3}
\left\{
\begin{array}{l}
\omega_T=\varepsilon|_T, \\
\omega_L-F_1'\omega_L^\sharp\lambda=\eta|_L+F_1'\omega_T^\sharp\tau.
\end{array}
\right.
\end{equation}

Let $\HHH$ be the abelian group of 
$\sigma$-linear homomorphisms $L\to\Fa P$. 
We claim that the endomorphism $U$ of $\HHH$ given by
$U(\omega_L)=F_1'\omega_L^\sharp\lambda$ is elementwise
nilpotent, which implies that $1-U$ is bijective, and
\eqref{EqCond3} has a unique solution.
The endomorphism $F_1'$ of $\Fa P$ is elementwise
nilpotent because $F_1'(ax)=\sigma_1(a)F'(x)$ and because
$\sigma_1$ is elementwise nilpotent on $\Fa$ by assumption.
Since $L$ is finitely generated it follows that
$U$ is elementwise nilpotent as desired.
\end{proof}

\begin{Remark}
The same argument applies if instead of $\sigma_1$ being element\-wise
nilpotent one demands that $\lambda$ is (topologically)
nilpotent, which is the original situation in 
\cite[Theorem 44]{Zink-Disp}; see section \ref{Se-disp}.
\end{Remark}

%----------------------------------------------------------------
\section{Abstract deformation theory}
\label{Se-abs-def}

\begin{Defn}
The Hodge filtration of a window $\PPP$
is the submodule 
$$
Q/IP\subseteq P/IP.
$$
\end{Defn}

\begin{Lemma}
\label{Le-Hodge}
Let $\alpha:\FFF\to\FFF'$ be a strict homomorphism
of frames with $S=S'$. Then $R\to R'$ is surjective
and we have $I\subseteq I'$.
Windows $\PPP$ over $\FFF$ are equivalent to pairs 
consisting of a window $\PPP'$ over $\FFF'$ and a lift 
of its Hodge filtration to a direct summand 
$V\subseteq P'/IP'$.
\end{Lemma}

\begin{proof}
The equivalence is given by the functor
$\PPP\mapsto(\alpha_*\PPP,Q/IP)$, 
which is easily seen to be fully faithful. 
We show that it is essentially surjective.
Let a window $\PPP'$ over $\FFF'$ and a lift if 
its Hodge filtration $V\subseteq P'/IP'$ 
be given and let $Q\subset P'$ be the 
inverse image of $V$. We have to show that
$\PPP=(P',Q,F',F_1'|_Q)$ is a window over $\FFF$.
First we need a normal decomposition for $\PPP$;
this is a decomposition $P'=L\oplus T$ such that $V=L/IL$.
Since $\PPP'$ has a normal decomposition, $\PPP$ 
has one too for at least one choice of $V$.
By modifying the isomorphism $P'\cong L\oplus T$
with an automorphism $\left(\begin{smallmatrix}
1 & 0\\ u & 1\end{smallmatrix}\right)$ of $L\oplus T$ 
for some homomorphism $u:L\to I'T$ one reaches 
every lift of the Hodge filtration.
It remains to show that $F_1'(Q)$ generates $P'$. 
In terms of a normal decomposition $P'=L\oplus T$ 
for $\PPP$ this means that $F_1'+F':L\oplus T\to P'$ 
is a $\sigma$-linear isomorphism, which holds because 
$\PPP'$ is a window.
\end{proof}

Assume that a strict homomorphism of frames $\alpha:\FFF\to\FFF'$
is given such that $S\to S'$ is surjective with kernel $\Fa$,
and $I'=IS'$.
We want to factor $\alpha$ into strict homomorphisms
\begin{equation}
\label{Eq-factor}
(S,I,R,\sigma,\sigma_1)\xrightarrow{\alpha_1}
(S,I'',R',\sigma,\sigma_1'')\xrightarrow{\alpha_2}
(S',I',R',\sigma',\sigma_1')
\end{equation}
such that $\alpha_2$ satisfies the hypotheses of
Theorem \ref{Th-crys}. 

Necessarily $I''=I+\Fa$. 
The main point is to define $\sigma_1'':I''\to S$, which is 
equivalent to defining a $\sigma$-linear homomorphism 
$\sigma_1'':\Fa\to\Fa$ that extends the restriction 
of $\sigma_1$ to $I\cap\Fa$ and satisfies 
the hypotheses of Theorem \ref{Th-crys}.

If this is achieved, Theorem \ref{Th-crys} and
Lemma \ref{Le-Hodge} show that windows over
$\FFF$ are equivalent to windows $\PPP'$ over $\FFF'$
plus a lift of the Hodge filtration
to a direct summand of $P/IP$, where
$\PPP''=(P,Q'',F,F_1'')$ is the unique lift 
of $\PPP'$ under $\alpha_2$.

%----------------------------------------------------------------
\section{Dieudonn\'e frames}

Let $R$ be a noetherian complete local ring with 
maximal ideal $\Fm$ and with perfect residue field 
$k$ of characteristic $p$. 
If $p=2$ we assume that $pR=0$. 
There is a unique subring $\WW(R)\subset W(R)$ 
stable under its Frobenius $f$ 
such that the projection $\WW(R)\to W(k)$
is surjective with kernel $\hat W(\Fm)$,
the ideal of all Witt vectors in $W(\Fm)$ 
whose coefficients converge to zero $\Fm$-adically,
and $\WW(R)$ is also stable under the Verschiebung $v$;
see \cite[Lemma 2]{Zink-DDisp}.
Let $\II_R$ be the kernel of the
projection to the first component $\WW(R)\to R$.
Then $v:\WW(R)\to\II_R$ is bijective.

\begin{Defn}
The Dieudonn\'e frame associated to $R$ is
$$
\DDD_R=(\WW(R),\II_R,R,f,f_1)
$$
with $f_1=v^{-1}$.
\end{Defn}

Here $\theta=p$.
Windows over $\DDD_R$ are Dieudonn\'e displays
over $R$ in the sense of \cite{Zink-DDisp}.
We note that $\WW(R)$ is a local ring, which
guarantees the existence of normal decompositions;
see Remark \ref{Re-lift}. 
The inclusion $\WW(R)\to W(R)$ is a homomorphism
of frames $\DDD_R\to\WWW_R$.
A local ring homomorphism $R\to R'$ induces a 
frame homomorphism $\DDD_R\to\DDD_{R'}$.

Assume that $R'=R/\Fb$ for an ideal $\Fb$ 
equipped with elementwise nilpotent divided powers.
Then $\WW(R)\to\WW(R')$ is surjective with 
kernel $\hat W(\Fb)=W(\Fb)\cap\hat W(\Fm)$.
In this situation, a factorisation \eqref{Eq-factor} 
of the homomorphism $\DDD_R\to\DDD_{R'}$ can be defined
as follows.

Let $\Fb^{<\infty>}$ be the $\WW(R)$-module of all
sequences $[b_0,b_1,\ldots]$ with elements $b_i\in\Fb$
that converge to zero $\Fm$-adically, on which
$x\in\WW(R)$ acts by 
$[b_0,b_1,\ldots]\mapsto [w_0(x)b_0,w_1(x)b_1,\ldots]$.
The divided Witt polynomials define an isomorphism
of $\WW(R)$-modules
$$
\log:\hat W(\Fb)\cong\Fb^{<\infty>}.
$$

Let $\tilde\II=\II_R+\hat W(\Fb)$.
In logarithmic coordinates, the restriction 
of $f_1$ to $\II_R\cap\hat W(\Fb)$ is given by
$$
f_1[0,b_1,b_2,\ldots]=[b_1,b_2,\ldots].
$$
Thus $f_1:\II_R\to\WW(R)$ extends uniquely 
to an $f$-linear homomorphism 
$$
\tilde f_1:\tilde\II\to\WW(R)
$$
with $\tilde f_1[b_0,b_1,\ldots]=[b_1,b_2,\ldots]$
on $\hat W(\Fb)$, and we obtain a factorisation
\begin{equation}
\label{Eq-Dieu-factor}
\DDD_R\xrightarrow{\alpha_1}
\DDD_{R/R'}=(\WW(R),\tilde\II,R',f,\tilde f_1)
\xrightarrow{\alpha_2}\DDD_{R'}.
\end{equation}

\begin{Prop}
\label{Pr-crys-Dieu}
The homomorphism $\alpha_2$ is crystalline.
\end{Prop}

This is a reformulation of \cite[Theorem 3]{Zink-DDisp}
if $\Fm$ is nilpotent, and the general case is an easy
consequence. As explained in section \ref{Se-abs-def},
it follows that deformations of Dieudonn\'e displays
from $R'$ to $R$ are classified by lifts of the Hodge
filtration; this is \cite[Theorem 4]{Zink-DDisp}.

\begin{proof}[Proof of Proposition \ref{Pr-crys-Dieu}]
When $\Fm$ is nilpotent, $\alpha_2$ satisfies the
hypotheses of Theorem \ref{Th-crys};
the required filtration of $\Fa=\hat W(\Fb)$ 
is $\Fa_i=p^i\Fa$. In general, these hypotheses
are not fulfilled
because $f_1:\Fa\to\Fa$ is only topologically nilpotent.
However, one can find a sequence of ideals
$R\supset I_1\supset I_2\ldots$ which define the
$\Fm$-adic topology such that each $\Fb\cap I_n$ is
is stable under the divided powers of $\Fb$. Indeed,
for each $n$ there is an $l$ with 
$\Fm^l\cap\Fb\subseteq\Fm^n\Fb$; for 
$I_n=\Fm^n\Fb+\Fm^l$ we have $\Fb\cap I_n=\Fm^n\Fb$.
The proposition holds for each $R/I_n$ in place of $R$,
and the general case follows by passing to the 
projective limit, using Lemma \ref{Le-lim}.
\end{proof}

%----------------------------------------------------------------
\section{$\varkappa$-frames}
\label{Se-kappa}

The results in this section are essentially due to Th.~Zink.

\begin{Defn}
A $\varkappa$-frame is a frame 
$\FFF=(S,I,R,\sigma,\sigma_1)$ such that
\begin{enumerate}
\setcounter{enumi}{\value{frame-count}}
\renewcommand{\theenumi}{\roman{enumi}}
\item \label{Def-kappa-tors-S}
$S$ has no $p$-torsion,
\item \label{Def-kappa-tors-W}
$W(R)$ has no $p$-torsion,
\item \label{Def-kappa-sigma-theta}
$\sigma(\theta)-\theta^p=p\cdot\text{unit}$ in $S$.
\end{enumerate}
\end{Defn}

\begin{Remarks}
If \eqref{Def-frame-Frob} and \eqref{Def-kappa-tors-S} hold
then we have a (non-additive) map 
$$
\tau:S\to S,\qquad\tau(x)=\frac{\sigma(x)-x^p}p,
$$
and \eqref{Def-kappa-sigma-theta} 
says that $\tau(\theta)$ is a unit.
Condition \eqref{Def-kappa-tors-W} is satisfied if and only if 
the nilradical $\NNN(R)$ 
has no $p$-torsion, for example if $R$ is reduced, 
or flat over $\ZZ_{(p)}$.
\end{Remarks}

\begin{Prop}
\label{Pr-kappa}
To each $\varkappa$-frame $\FFF$ one can associate
a unit $u$ of $W(R)$ and a $u$-homomorphism of frames 
$\varkappa:\FFF\to\WWW_R$ lying over $\id_R$.
The construction is functorial in $\FFF$.
\end{Prop}

\begin{proof}
Condition \eqref{Def-kappa-tors-S} implies that there is 
a well-defined homomorphism $\delta:S\to W(S)$ 
with $w_n\delta=\sigma^n$; 
see \cite[IX.1, proposition 2]{Bour-Comm-Alg}. 
We have $f\delta=\delta\sigma$.
Let $\varkappa$ be the composite ring homomorphism
$$
\varkappa:S\xrightarrow{\delta}W(S)\to W(R).
$$
Then $f\varkappa=\varkappa\sigma$
and $\varkappa(I)\subseteq I_R$.
Clearly $\varkappa$ is functorial in $\FFF$.

To define $u$ we write $1=\sum y_i\sigma_1(x_i)$
in $S$ with $x_i\in I$ and $y_i\in S$. 
This is possible by \eqref{Def-frame-surj}.
We recall that $\theta=\sum y_i\sigma(x_i)$;
see the proof of Lemma \ref{Le-theta}.
Let $u=\sum\varkappa(y_i)f_1\varkappa(x_i)$.
Then $pu=\varkappa(\theta)$ because $pf_1=f$.
We claim that $f_1\varkappa=u\cdot\varkappa\sigma_1$.
By \eqref{Def-kappa-tors-W} this is equivalent
to $p\cdot f_1\varkappa=pu\cdot\varkappa\sigma_1$,
which is easily checked as $pf_1=f$ and 
$\theta\sigma_1=\sigma$.

It remains to show that $u$ is a unit in $W(R)$.
Let $pu=\varkappa(\theta)=(a_0,a_1,\ldots)$ 
as a Witt vector.
By Lemma \ref{Le-unit} below, $u$ is a unit 
if and only if $a_1$ is a unit in $R$.
But $\delta(\theta)=(\theta,\tau(\theta),\ldots)$
because $w_*$ applied to both sides gives 
$(\theta,\sigma(\theta),\ldots)$; here `$\ldots$' means
`not specified'.
Hence $a_1$ is a unit by \eqref{Def-kappa-sigma-theta}.

Finally, $u$ is functorial in $\FFF$ by its uniqueness, 
see Remark \ref{Re-u-unique}.
\end{proof}

\begin{Lemma}
\label{Le-unit}
Let $R$ be a ring with $p\in\Rad(R)$ and let $u\in W(R)$.
For an integer $r\geq 0$ let $p^ru=(a_0,a_1,a_2,\ldots)$.
Then $u$ is a unit in $W(R)$ if and only if $a_r$ is a unit in $R$.
\end{Lemma}

\begin{proof}
Let $r=0$. It suffices to show that an element
$\bar u\in W_{n+1}(R)$ that maps to $1$ in $W_n(R)$
is a unit. If $\bar u=1+v^n(x)$ with $x\in R$
then $\bar u^{-1}=1+v^n(y)$ where $y\in R$ is
determined by $x+y+pxy=0$, which has a solution
as $p\in\Rad(R)$. For general $r$, by the case $r=0$ 
we may replace $R$ by $R/pR$. Then we have 
$p(b_0,b_1,\ldots)=(0,b_0^p,b_1^p,\ldots)$ in $W(R)$,
which reduces the assertion to the case $r=0$.
\end{proof}

\begin{Cor}
\label{Co-kappa}
Let $\FFF$ be a $\varkappa$-frame with
$S=W(k)[[x_1,\ldots,x_r]]$ for a perfect field $k$ 
of odd characteristic $p$. Assume that $\sigma$ extends
the Frobenius automorphism of $W(k)$ by $\sigma(x_i)=x_i^p$. 
Then $u$ is a unit in $\WW(R)$ and $\varkappa$ induces
a $u$-homomorphism of frames $\varkappa:\FFF\to\DDD_R$.
\end{Cor}

\begin{proof}
We claim that $\delta(S)$ lies in $\WW(S)$.
Indeed, $\delta(x_i)=[x_i]$ because $w_n$ applied to
both sides gives $x_i^{p^n}$. 
Thus $\delta(x^e)=[x^e]\in\WW(S)$ 
for any multi-exponent $e=(e_1,\ldots e_r)$. 
Since $\WW(S)=\varprojlim\WW(S/\Fm^n)$ and since for each $n$
all but finitely many $x^e$ lie in $\Fm^n$ the claim follows.
Hence the image of $\varkappa:S\to W(R)$ is contained in $\WW(R)$. 
By construction the element $u$ lies in $\WW(R)$; 
it is invertible in $\WW(R)$ because the inclusion 
$\WW(R)\to W(R)$ is a local homomorphism of local rings.
\end{proof}

%----------------------------------------------------------------
\section{The main frame}
\label{Se-main-frame}

Let $R$ be a complete regular local ring with perfect
residue field $k$ of characteristic $p\geq 3$. 
We choose a continuous ring homomorphism
$$
\FS=W(k)[[x_1,\ldots,x_r]]\xrightarrow{\pi} R
$$
such that $x_1,\ldots,x_r$ map to a regular
system of parameters of $R$.
As the graded ring of $R$ is isomorphic to
$k[x_1,\ldots,x_r]$, one can find a power series 
$E_0\in\FS$ with constant term zero
such that $\pi(E_0)=-p$. Let $E=E_0+p$
and $I=E\FS$. Then $R=\FS/I$.
Let $\sigma:\FS\to\FS$ be the continuous ring endomorphism
that extends the Frobenius automorphism of $W(k)$ by 
$\sigma(x_i)=x_i^p$. We have a frame
$$
\BBB=(\FS,I,R,\sigma,\sigma_1)
$$
where $\sigma_1(Ey)=\sigma(y)$ for $y\in\FS$.

\begin{Lemma}
\label{Le-B-kappa}
The frame $\BBB$ is a $\varkappa$-frame.
\end{Lemma}

\begin{proof}
Let $\theta\in\FS$ be the element given by Lemma \ref{Le-theta}.
The only condition to be checked is that $\tau(\theta)$
is a unit in $\FS$.  Let $E_0'=\sigma(E_0)$.
As $\sigma_1(E)=1$ we have $\theta=\sigma(E)=E_0'+p$. Hence 
$$
\tau(\theta)=\frac{\sigma(E_0')+p-(E_0'+p)^p}p
\equiv 1+\tau(E_0')\mod p.
$$
Since the constant term of $E_0$ is zero, the same is true
for $\tau(E_0')$, which implies that $\tau(\theta)$ is a unit
as required.
\end{proof}

By Proposition \ref{Pr-kappa} and Corollary \ref{Co-kappa} 
we get a ring homomorphism $\varkappa:\FS\to\WW(R)$, 
which is a $u$-homomorphism of frames
$$
\varkappa:\BBB\to\DDD_R.
$$
Here the unit $u\in\WW(R)$ is determined by $pu=\varkappa\sigma(E)$.

\begin{Thm}
\label{Th-kappa-crys}
The homomorphism $\varkappa$ is crystalline 
(Definition \ref{Def-crys}).
\end{Thm}

To prove this we consider the following auxiliary frames.
Let $J\subset\FS$ be the ideal $J=(x_1,\ldots,x_r)$.
For $a\in\NN$ let $\FS_a=\FS/J^a\FS$
and let $R_a=R/\Fm_R^a$. Then $R_a=\FS_a/E\FS_a$.
The element $E$ is not a zero divisor in $\FS_a$. 
There is a well-defined frame
$$
\BBB_a=(\FS_a,I_a,R_a,\sigma_a,\sigma_{1a})
$$
such that the projection $\FS\to\FS_a$ is a strict
homomorphism $\BBB\to\BBB_a$. Indeed, $\sigma$
induces an endomorphism $\sigma_a$ of $\FS_a$ 
because $\sigma(J)\subseteq J$, and for
$y\in\FS_a$ one can define $\sigma_{1a}(Ey)=\sigma_a(y)$.

For simplicity, the image of $u$ in $\WW(R_a)$
is denoted by $u$ as well. The $u$-homomorphism
$\varkappa$ induces a $u$-homomorphism
$$
\varkappa_a:\BBB_a\to\DDD_{R_a}
$$
because for $e\in\NN^r$ we have $\varkappa(x^e)=[x^e]$, 
which maps to zero in $\WW(R_a)$ 
when $e_1+\ldots+e_r\geq a$. 
We note that $\BBB_a$ is again a $\varkappa$-frame, 
so the existence of $\varkappa_a$ can also be viewed 
as a consequence of Proposition \ref{Pr-kappa}.

\begin{Thm}
\label{Th-kappa-a-crys}
For each $a\in\NN$ the homomorphism $\varkappa_a$
is crystalline.
\end{Thm}

This includes Theorem \ref{Th-kappa-crys} if one allows
$a=\infty$ and writes $\BBB=\BBB_\infty$ etc.
To prepare for the proof, for each $a\in\NN$
we want to construct the following commutative diagram 
of frames where vertical arrows are $u$-homomorphisms 
and where horizontal arrows are strict.

\begin{equation}
\label{Diag}
\xymatrix@!C@M+0.2em{
\BBB_{a+1} \ar[r] \ar[d]^{\varkappa_{a+1}} &
\tBBB_{a+1} \ar[r]^-{\pi} \ar[d]^{\tvarkappa_{a+1}} &
\BBB_a \ar[d]^{\varkappa_a} \\
\DDD_{R_{a+1}} \ar[r] & 
\DDD_{R_{a+1}/R_a} \ar[r]^-{\pi'} & \DDD_{R_a} \hspace{-0.5em}
}
\end{equation}

The upper line is a factorisation \eqref{Eq-factor} 
of the projection $\BBB_{a+1}\to\BBB_a$. 
This means that the frame $\tBBB_{a+1}$ 
necessarily takes the form
$$
\tBBB_{a+1}=(\FS_{a+1},\tI_{a+1},R_a,\sigma_{a+1},\tsigma_{1(a+1)})
$$
with $\tI_{a+1}=E\FS_{a+1}+J^a/J^{a+1}$.
We define $\tsigma_{1(a+1)}:\tI_{a+1}\to\FS_{a+1}$ to be 
the extension of $\sigma_{1(a+1)}:E\FS_{a+1}\to\FS_{a+1}$ 
by zero on $J^a/J^{a+1}$. This is well-defined because 
$$
E\FS_{a+1}\cap J^a/J^{a+1}=E(J^a/J^{a+1})
$$
and because for $x\in J^a/J^{a+1}$ we have
$\sigma_{1(a+1)}(Ex)=\sigma_{a+1}(x)$,
which is zero as $\sigma(J^a)\subseteq J^{ap}$.

The lower line of \eqref{Diag} is the factorisation 
\eqref{Eq-Dieu-factor} with respect to the trivial
divided powers on the kernel $\Fm_R^a/\Fm_R^{a+1}$.

In order that the diagram commutes it is necessary
and sufficient that $\tvarkappa_{a+1}$ is given by
the ring homomorphism $\varkappa_{a+1}$.

It remains to show that $\tvarkappa_{a+1}$ is a 
$u$-homomorphism of frames. 
The only non-trivial condition is that
$\tf_1\varkappa_{a+1}=u\cdot\varkappa_{a+1}\tsigma_{1(a+1)}$ 
on $\tI_{a+1}$. This relation holds on $E\FS_{a+1}$ because 
$\varkappa_{a+1}$ is a $u$-homomorphism of frames.
On $J^a/J^{a+1}$ we have $\varkappa_{a+1}\tsigma_{1(a+1)}=0$ 
by definition. For $y\in\FS_{a+1}$ and for $e\in\NN^r$ with 
$e_1+\ldots+e_r=a$ we compute
$$
\tf_1(\varkappa_{a+1}(x^ey))
=\tf_1([x^e]\varkappa_{a+1}(y))
=\tf_1([x^e])f(\varkappa_{a+1}(y))=0
$$
because $\log([x^e])=\left<x^e,0,0,\ldots\right>$.
As these $x^e$ generate $J^a$, the required relation on
$J^a/J^{a+1}$ follows, and the diagram is constructed.

\begin{proof}[Proof of Theorem \ref{Th-kappa-a-crys}]
We use induction on $a$. The homomorphism $\varkappa_1$ 
is crystalline because it is bijective.
Assume that $\varkappa_a$ is crystalline for some
$a\in\NN$ and consider the diagram \eqref{Diag}.
The homomorphism $\pi'$ is crystalline by 
Proposition \ref{Pr-crys-Dieu}, while $\pi$ is crystalline 
by Theorem \ref{Th-crys}; the required filtration of
$J^a/J^{a+1}$ is trivial. Hence $\tvarkappa_{a+1}$ is crystalline.
By Lemma \ref{Le-Hodge} it follows that $\varkappa_{a+1}$
is crystalline too.
\end{proof}

\begin{proof}[Proof of Theorem \ref{Th-kappa-crys}] 
Use Theorem \ref{Th-kappa-a-crys} and Lemma \ref{Le-lim}. 
\end{proof}

%----------------------------------------------------------------
\section{Classification of group schemes}

The following consequences of Theorem \ref{Th-kappa-crys}
are analogous to \cite{Vasiu-Zink-Classif}.
Let $\BBB=(\FS,I,R,\sigma,\sigma_1)$ be the frame defined 
in section \ref{Se-main-frame}.

\begin{Defn}
A Breuil window relative to $\FS\to R$ is a pair
$(Q,\phi)$ where $Q$ is a free $\FS$-module of finite
rank and where $\phi:Q\to Q^{(\sigma)}$ is an $\FS$-linear 
homomorphism with cokernel annihilated by $E$.
\end{Defn}

\begin{Lemma}
\label{Le-Br-win}
Breuil windows relative to $\FS\to R$ are equivalent
to $\BBB$-windows in the sense of Definition \ref{Def-win}.
\end{Lemma}

\begin{proof}
This is similar to \cite[Lemma 1]{Vasiu-Zink-Classif}.
For a window $(P,Q,F,F_1)$ over $\BBB$ the module 
$Q$ is free over $\FS$ because $I=E\FS$ is free.
Hence $F_1^\sharp:Q^{(\sigma)}\to P$ is bijective, and
we can define a Breuil window $(Q,\phi)$ where
$\phi$ is the inclusion $Q\to P$
composed with the inverse of $F_1^\sharp$. 
Conversely, if $(Q,\phi)$ is a Breuil window, 
$\Coker(\phi)$ is a free $R$-module.
Indeed, $\phi$ is injective because it becomes 
bijective over $\FS[E^{-1}]$, so $\Coker(\phi)$
has projective dimension one over $\FS$, which 
implies that it is free over $R$ by using depth.
Thus one can define a window over $\BBB$ as 
follows: $P=Q^{(\sigma)}$, the inclusion $Q\to P$
is $\phi$, $F_1:Q\to Q^{(\sigma)}$ is the homomorphism 
$x\mapsto 1\otimes x$, and $F(x)=F_1(Ex)$.
The two constructions are mutually inverse.
\end{proof}

By \cite{Zink-DDisp}, $p$-divisible groups over $R$ 
are equivalent to Dieudonn\'e displays over $R$.
Together with Theorem \ref{Th-kappa-crys} and 
Lemma \ref{Le-Br-win} this implies:

\begin{Cor}
\label{Co-windows}
The category of $p$-divisible groups over $R$ is equivalent 
to the category of Breuil windows relative to $\FS\to R$.
\qed
\end{Cor}

Let us use the following abbreviation:
An {\em admissible torsion module} is a finitely generated 
$\FS$-module annihilated by a power of $p$ and of 
projective dimension at most one.

\begin{Defn}
A Breuil module relative to $\FS\to R$ is a triple 
$(M,\varphi,\psi)$ where $M$ is an admissible torsion module
together with $\FS$-linear homomorphisms 
$\varphi:M\to M^{(\sigma)}$ and $\psi:M^{(\sigma)}\to M$ 
such that $\varphi\psi=E$ and $\psi\varphi=E$.
\end{Defn}

When $R$ has characteristic zero, each of the maps $\varphi$ 
and $\psi$ determines the other one; see 
Lemma \ref{Le-Br-Mod} below.

\begin{Thm}
\label{Th-modules}
The category of finite flat group schemes 
over $R$ annihilated by a power of $p$ is equivalent to the
category of Breuil modules relative to $\FS\to R$.
\end{Thm}

\begin{proof}
The assertion follows from Corollary \ref{Co-windows} by the 
arguments of \cite{Kisin} or \cite{Vasiu-Zink-Classif},
which we recall briefly.
A homomorphism of Breuil windows $(Q',\phi')\to(Q,\phi)$ 
is called an isogeny if it becomes
bijective over $\FS[1/p]$. Then its cokernel
is naturally a Breuil module; the required $\psi$ is
induced by the homomorphism $E\phi^{-1}:Q^{(\sigma)}\to Q$.
Using that the equivalence between $p$-divisible groups
and Breuil windows preserves isogenies and short exact sequences,
Theorem \ref{Th-modules} is a formal consequence
of the following two facts.

\smallskip
(a) 
Each finite flat group scheme over $R$ of $p$-power order 
is the kernel of an isogeny of $p$-divisible groups. 
See \cite[Th\'eor\`eme 3.1.1]{BBM}.

\smallskip
(b)
Each Breuil module relative to $\FS\to R$ is the cokernel 
of an isogeny of Breuil windows. This is analogous to 
\cite[Proposition 2]{Vasiu-Zink-Classif}.

\smallskip
Let us recall the argument for (b).
If $(M,\varphi,\psi)$ is a Breuil module,
one can find free $\FS$-modules $P$ and $Q$ together
with surjections $Q\to M$ and $P\to M^{(\sigma)}$
and homomorphisms $\tvarphi:Q\to P$ and
$\tpsi:P\to Q$ which lift $\varphi$ and $\psi$
such that $\tvarphi\tpsi=E$ and $\tpsi\tvarphi=E$.
Next one chooses an isomorphism $\alpha:P\cong Q^{(\sigma)}$
compatible with the given projections of both sides to 
$M^{(\sigma)}$. Let $\phi=\alpha\tvarphi$.
Then $(M,\varphi,\psi)$ is the cokernel of the isogeny
$(Q',\phi')\to(Q,\phi)$, where $Q'$ is the kernel
of $Q\to M$ and $\phi'$ is the restriction of $\phi$.
\end{proof}

\begin{Lemma}
\label{Le-Br-Mod}
If $R$ has characteristic zero, the category of Breuil
modules relative to $\FS\to R$ is equivalent to the
category of pairs $(M,\varphi)$ where $M$ is an 
admissible torsion module and where 
$\varphi:M\to M^{(\sigma)}$ is an $\FS$-linear 
homomorphism with cokernel annihilated by $E$.
\end{Lemma}

\begin{proof}
Cf.\ \cite[Proposition 2]{Vasiu-Zink-Classif}.
For a non-zero admissible torsion module $M$ the 
set of zero divisors on $M$ is equal to $\Fp=p\FS$
because every associated prime of $M$ has height one
and contains $p$.
In particular, $M\to M_\Fp$ is injective.
The hypothesis of the lemma means that $E\not\in\Fp$.
For a given pair $(M,\varphi)$ as in the lemma this 
implies that $\varphi_\Fp:M_\Fp\to M_\Fp^{(\sigma)}$ 
is surjective, thus bijective because both sides
have the same finite length. It follows that 
$\varphi$ is injective, and $(M,\varphi)$ is
extended uniquely to a Breuil module by
$\psi(x)=\varphi^{-1}(Ex)$. 
\end{proof}

\subsection*{Duality}

The dual of a Breuil window $(Q,\phi)$ is the Breuil window 
$(Q,\phi)^t=(Q^\vee,\psi^\vee)$ where $Q^\vee=\Hom_\FS(Q,\FS)$
and where $\psi:Q^{(\sigma)}\to Q$ is the unique homomorphism 
with $\psi\phi=E$. Here we identify $(Q^{(\sigma)})^\vee$
and $(Q^\vee)^{(\sigma)}$. 
For a $p$-divisible group $G$ over $R$ let 
$G^\vee$ be the Cartier dual of $G$ and let 
$\MM(G)$ be the Breuil window associated with $G$ 
by the equivalence of Corollary \ref{Co-windows}.

\begin{Prop}
\label{Pr-dual-windows}
There is a natural isomorphism $\MM(G^\vee)\cong\MM(G)^t$.
\end{Prop}

\begin{proof}
The equivalence between $p$-divisible groups over $R$ and
Dieudon\-n\'e displays over $R$ is compatible with duality
by \cite[Theorem 3.4]{Lau-Dual}. It is easy to see that the
equivalence of Lemma \ref{Le-Br-win} and the functor
$\varkappa_*$ preserve duality as well. The proposition follows.
\end{proof}

The dual of a Breuil module $\MM=(M,\varphi,\psi)$ is 
the Breuil module $\MM^t=(M^\star,\psi^\star,\varphi^\star)$ 
where $M^\star=\Ext^1_\FS(M,\FS)$. Here we identify 
$(M^{(\sigma)})^\star$ and $(M^\star)^{(\sigma)}$
using that $(\;)^{(\sigma)}$ preserves projective
resolutions as $\sigma$ is flat.
For a finite flat group scheme $H$ over $R$ of $p$-power 
order let $H^\vee$ be the Cartier dual of $H$ and let 
$\MM(H)$ be the Breuil module associated with $H$ by
the equivalence of Theorem \ref{Th-modules}.

\begin{Prop}
\label{Pr-dual-modules}
There is a natural isomorphism $\MM(H^\vee)\cong\MM(H)^t$.
\end{Prop}

\begin{proof}
Choose an isogeny of $p$-divisible groups $G_1\to G_2$
with kernel $H$. Then $\MM(H)$ is the cokernel of 
$\MM(G_1)\to\MM(G_2)$, which implies that $\MM(H)^t$ 
is the cokernel of $\MM(G_2)^t\to\MM(G_1)^t$.
On the other hand, $H^\vee$ is the kernel of 
$G_2^\vee\to G_1^\vee$, so $\MM(H^\vee)$ is the cokernel
of $\MM(G_2^\vee)\to\MM(G_1^\vee)$. 
Proposition \ref{Pr-dual-windows} applied to $G_1$ and $G_2$ 
gives an isomorphism $\beta:\MM(H^\vee)\cong\MM(H)^t$. 
One easily checks that $\beta$ is independent of the 
choice and functorial in $H$.
\end{proof}

%----------------------------------------------------------------
\section{Other lifts of Frobenius}

One may ask how much freedom we have in the
choice of $\sigma$ for the frame $\BBB$. 
Let $J=(x_1,\ldots,x_r)$. To begin with, let $\sigma:\FS\to\FS$
be an arbitrary ring endomorphism such that $\sigma(J)\subset J$ 
and $\sigma(a)\equiv a^p$ mod $p\FS$ for $a\in\FS$. 
As in Section \ref{Se-main-frame} we consider the frame
$$
\BBB=(\FS,I,R,\sigma,\sigma_1)
$$
with $\sigma_1(Ey)=\sigma(y)$. Again this is
a $\varkappa$-frame because the proof of 
Lemma \ref{Le-B-kappa} uses only that $\sigma$
preserves $J$, so Proposition \ref{Pr-kappa}
gives a homomorphism of frames
$$
\varkappa:\BBB\to\WWW_R.
$$
By the assumptions on $\sigma$ we have 
$\sigma(J)\subseteq J^p+pJ$, which implies that the
endomorphism $\sigma:J/J^2\to J/J^2$ is divisible by $p$.

\begin{Prop}
\label{Pr-cond-sigma}
The image of $\varkappa:\FS\to W(R)$ lies in $\WW(R)$ 
if and only if the endomorphism $\sigma/p$ of $J/J^2$ 
is nilpotent modulo $p$.
\end{Prop}

We have a non-additive map $\tau:J\to J$
given by $\tau(x)=(\sigma(x)-x^p)/p$.
Let $\Fm\subset\FS$ be the maximal ideal.
We write $gr_n(J)=\Fm^nJ/\Fm^{n+1}J$.

\begin{Lemma}
\label{Le-cond-sigma}
For $n\geq 0$ the map $\tau$ preserves $\Fm^nJ$ and induces a
$\sigma$-linear endomorphism of $k$-modules
$gr_n(\tau):gr_n(J)\to gr_n(J)$.
We have $gr_0(\tau)=\sigma/p$ as an endomorphism of
$gr_0(J)=J/J^2+pJ$.
There is a commutative diagram of the following type 
with $\pi i=\id$.
$$
\xymatrix@M+0.2em{
gr_n(J) \ar[r]^{gr_n(\tau)} \ar[d]_{\pi} & gr_n(J) \\ 
gr_0(J) \ar[r]^{gr_0(\tau)} & gr_0(J) \ar[u]_i
}
$$
\end{Lemma}

\begin{proof}
Let $J'=p^{-1}\Fm J$ as an $\FS$-submodule of $J\otimes\QQ$.
Then $J\subset J'$, and 
$gr_n(J)$ is a submodule of $gr_n(J')=\Fm^nJ'/\Fm^{n+1}J'$.
The composition $J\xrightarrow{\tau}J\subset J'$ 
can be written as
$\tau=\sigma/p-\varphi/p$, where $\varphi(x)=x^p$.
One checks that $\varphi/p:\Fm^nJ\to\Fm^{n+1}J'$
(which requires $p\geq 3$ when $n=0$) and that
$\sigma/p:\Fm^nJ\to\Fm^nJ'$.  
Hence $\sigma/p$ and $\tau$ induce 
the same map $\Fm^nJ\to gr_n(J')$. This map
is $\sigma$-linear and zero on $\Fm^{n+1}J$
because this holds for $\sigma/p$,
and its image lies in $gr_n(J)$ because 
this is true for $\tau$.

We define $i:gr_0(J)\to gr_n(J)$ by
$x\mapsto p^n x$. For $n\geq 1$ let $K_n$ be the image 
of $\Fm^{n-1}J^2\to gr_n(J)$. Then $i$ maps $gr_0(J)$
bijectively onto $gr_n(J)/K_n$, so there is a unique
homomorphism $\pi:gr_n(J)\to gr_0(J)$ with kernel $K_n$ 
such that $\pi i=\id$. Since $i$ commutes with $gr(\tau)$,
in order that the diagram commutes it suffices that $gr_n(\tau)$ 
vanishes on $K_n$. We have $\sigma(J)\subseteq\Fm J$, which 
implies that $(\sigma/p)(\Fm^{n-1}J^2)\subseteq\Fm^{n+1}J'$,
and the assertion follows.
\end{proof}

\begin{proof}[Proof of Proposition \ref{Pr-cond-sigma}]
We recall that $\varkappa=\pi\delta$, where 
$\delta:\FS\to W(\FS)$ is defined by 
$w_n\delta=\sigma^n$ for $n\geq 0$, 
and where $\pi:W(\FS)\to W(R)$ is the obvious projection.
For $x\in J$ and $n\geq 1$ let 
$$
\tau_n(x)=(\sigma(x)^{p^{n-1}}-x^{p^{n}})/p^n,
$$
thus $\tau_1=\tau$. It is easy to see that 
$$
\tau_{n+1}(x)\in J\cdot\tau_n(x),
$$
in particular we have $\tau_n:J\to J^n$.
If $\delta(x)=(y_0,y_1,\ldots)$,
the coefficients $y_n$ are determined by $y_0=x$ and
$w_n(y)=\sigma w_{n-1}(y)$ for $n\geq 1$, 
which translates into the equations
$$
y_n=\tau_n(y_0)+\tau_{n-1}(y_1)+\ldots+\tau_1(y_{n-1}).
$$

Assume now that $\sigma/p$ is nilpotent on $J/J^2$ modulo $p$.
By Lemma \ref{Le-cond-sigma} this implies that $gr_n(\tau)$ 
is nilpotent for every $n\geq 0$. We have to show that
for $x\in J$ the element $\delta(x)$ lies is $\WW(\FS)$, 
which means that the above sequence $(y_n)$ converges to zero. 
Assume that for some $N\geq 0$ we have $y_n\in\Fm^NJ$
for all but finitely many $n$.
The last two displayed equations give that 
$y_n-\tau(y_{n-1})\in\Fm^{N+1}J$ for all but finitely many $n$. 
As $gr_N(\tau)$ is nilpotent it follows that
$y_n\in\Fm^{N+1}J$ for all but finitely many $n$. 
Thus $\delta(x)\in\WW(\FS)$ and
in particular $\varkappa(x)\in\WW(R)$.

Conversely, if $\sigma/p$ is not nilpotent on $J/J^2$ modulo $p$,
then $gr_0(\tau)$ is not nilpotent by Lemma \ref{Le-cond-sigma},
so there is an $x\in J$ such that $\tau^n(x)\not\in\Fm J$ for
all $n\geq 0$. For $\delta(x)=(y_0,y_1,\ldots)$ 
we have $y_n\equiv\tau^nx$ modulo
$\Fm J$. The projection $\FS\to R$ induces an
isomorphism $J/\Fm J\cong\Fm_R/\Fm_R^2$. It follows
that $\varkappa(x)$ lies in $W(\Fm_R)$ but not in
$\hat W(\Fm_R)$, thus $\varkappa(x)\not\in\WW(R)$.
\end{proof}

Now we assume that $\sigma/p$ is nilpotent on 
$J/J^2$ modulo $p$. Then we have a homomorphism of frames 
$$
\varkappa:\BBB\to\DDD_R.
$$
As earlier let $\BBB_a=(\FS_a,I_a,R_a,\sigma_a,\sigma_{1a})$
with $\FS_a=\FS/J^a$ and with $R_a=R/\Fm_R^a$. 
The proof of Lemma \ref{Le-B-kappa} shows that $\BBB_a$
is a $\varkappa$-frame. Since $\WW(R_a)$ is the
image of $\WW(R)$ in $W(R_a)$ we get a homomorphism
of frames compatible with $\varkappa$,
$$
\varkappa_a:\BBB_a\to\DDD_{R_a}.
$$

\begin{Thm}
\label{Th-kappa-crys-gen}
The homomorphisms $\varkappa$ and $\varkappa_a$ are crystalline.
\end{Thm}

\begin{proof}
The proof is similar to that of Theorems 
\ref{Th-kappa-crys} and \ref{Th-kappa-a-crys}.

First we repeat the construction of \eqref{Diag}.
The restriction of $\sigma_{1(a+1)}$ to 
$E(J^a/J^{a+1})=p(J^a/J^{a+1})$ is given by 
$\sigma_1=\sigma/p=\tau$, which need not be zero 
in general, but still $\sigma_1$ extends uniquely 
to $J^a/J^{a+1}$ by the formula $\sigma_1=\sigma/p$. 
In order that $\tvarkappa_{a+1}$ is a $u$-homomorphism 
we have to check that
$\tf_1\varkappa_{a+1}=u\cdot \varkappa_{a+1}\tsigma_{1(a+1)}$
on $J^a/J^{a+1}$. Here $u$ acts on $J^a/J^{a+1}$ as the
identity. By the proof of Proposition \ref{Pr-cond-sigma},
for $x\in J^a/J^{a+1}$ we have
$$
\delta(x)=(x,\tau(x),\tau^2(x),\ldots).
$$
Since $\tsigma_{1(a+1)}(x)=\tau(x)$ the required relation follows.

To complete the proof we have to show that 
$\pi:\tBBB_{a+1}\to\BBB_a$ is crystalline.
Now $\sigma/p$ is nilpotent modulo $p$ on $J^n/J^{n+1}$
for $n\geq 1$.
Indeed, for $n=1$ this is our assumption, and for $n\geq 2$
the endomorphism $\sigma/p$ of $J^n/J^{n+1}$ is divisible
by $p^{n-1}$ as $\sigma(J)\subseteq pJ+J^p$.
In order to apply Theorem \ref{Th-crys} we need another
sequence of auxiliary frames: 
For $c\in\NN$ let $\FS_{a+1,c}=\FS_{a+1}/p^cJ^{a}\FS_{a+1}$
and let $\tBBB_{a+1,c}=(\FS_{a+1,c},I_{a+1,c},R_{a},\ldots)$
be the obvious quotient frame of $\tBBB_{a+1}$.
Then $\BBB_a$ is isomorphic to $\tBBB_{a+1,0}$, 
and $\tBBB_{a+1}$ is the projective limit of 
$\tBBB_{a+1,c}$ for $c\to\infty$.
Theorem \ref{Th-crys} shows that each projection 
$\tBBB_{a+1,c+1}\to\tBBB_{a+1,c}$ is crystalline,
which implies that $\pi$ is crystalline by Lemma \ref{Le-lim}.
\end{proof}

If $\sigma/p$ is nilpotent on $J/J^2$ modulo $p$, then
Corollary \ref{Co-windows}, Theorem \ref{Th-modules} 
and the Duality Propositions \ref{Pr-dual-windows} 
and \ref{Pr-dual-modules} follow as before.

%----------------------------------------------------------------
\section{Nilpotent windows}
\label{Se-disp}

All results in this article have a nilpotent counterpart
where only connected $p$-divisible groups and nilpotent 
windows are considered; then $k$ need not be perfect
and $p$ need not be odd. The necessary modifications
are standard, but for completeness we work out the details. 

\subsection{Nilpotence condition}

Let $\FFF=(S,I,R,f,f_1)$ be a frame.
For an $\FFF$-window $\PPP$ there is
a unique homomorphism of $S$-modules
$$
V^\sharp:P^{(\sigma)}\to P
$$
with $V^\sharp(F_1(x))=1\otimes x$ for $x\in Q$.
In terms of a normal representation $\Psi:L\oplus T\to P$
of $\PPP$ we have $V^\sharp=(1\oplus\theta)(\Psi^\sharp)^{-1}$.
The composition 
$$
P\xrightarrow{V^\sharp}P^{(\sigma)}
\xrightarrow{(V^\sharp)^{(\sigma)}}P^{(\sigma^2)}
\to \ldots \to P^{(\sigma^n)}
$$
is denoted $(V^\sharp)^n$ for simplicity.
The nilpotence condition
depends on the choice of an ideal $J\subset S$ 
such that $\sigma(J)+I+\theta S\subseteq J$, which we
call an {\em ideal of definition} for $\FFF$.

\begin{Defn}
Let $J\subset S$ an ideal of definition for $\FFF$. 
An $\FFF$-window $\PPP$ is called nilpotent
(with respect to $J$) if $(V^\sharp)^n\equiv 0$
modulo $J$ for sufficiently large $n$.
\end{Defn}

\begin{Remark}
For an $\FFF$-window $\PPP$ we consider the composition
$$
\lambda:L\subseteq L\oplus T\xrightarrow
{(\Psi^{\sharp})^{-1}}L^{(\sigma)}\oplus T^{(\sigma)}\to L^{(\sigma)}.
$$
Then $\PPP$ is nilpotent if and only if $\lambda$ 
is nilpotent modulo $J$.
\end{Remark}

\subsection{Nil-crystalline homomorphisms}

If $\alpha:\FFF\to\FFF'$ is a homomorphism of frames
and $J\subset S$ and $J'\subset S'$ are ideals of
definition with $\alpha(J)\subseteq J'$, the
functor $\alpha_*$ preserves nilpotent windows. 
We call $\alpha$ nil-crystalline if it induces an equivalence
between nilpotent $\FFF$-windows and nilpotent $\FFF'$-windows. 
The following variant of Theorem \ref{Th-crys} formalises
\cite[Theorem 44]{Zink-Disp}.

\begin{Thm}
\label{Th-nil-crys}
Let $\alpha:\FFF\to\FFF'$ be a homomorphism of frames
that induces an isomorphism $R\cong R'$ and a surjection
$S\to S'$ with kernel $\Fa\subset S$.
We assume that there is a finite filtration 
$\Fa=\Fa_0\supseteq\ldots\supseteq\Fa_n=0$
such that $\sigma(\Fa_i)\subseteq\Fa_{i+1}$ and 
$\sigma_1(\Fa_i)\subseteq\Fa_i$. 
We assume that finitely generated projective $S'$-modules
lift to projective $S$-modules. If $J\subset S$
is an ideal of definition such that $J^n\Fa=0$
for large $n$, then $\alpha$ is nil-crystalline 
with respect to $J\subset S$ and $J'=J/\Fa\subset S'$.
\end{Thm}

\begin{proof}
The assumptions imply that $\Fa\subseteq I\subseteq J$,
in particular $J'$ is well-defined.
An $\FFF$-window $\PPP$ is nilpotent
if and only if $\alpha_*\PPP$ is nilpotent. 
Using this, the proof of Theorem \ref{Th-crys} applies 
with the following modification in the final paragraph.
We claim that the endomorphism $U$ of $\HHH$ is nilpotent,
which again implies that $1-U$ is bijective.
Since $\PPP$ is nilpotent, $\lambda$ is nilpotent 
modulo $J$, so $\lambda$ is nilpotent modulo $J^n$
for each $n\geq 1$ as $J$ is stable under $\sigma$.
Since $J^n\Fa=0$ by assumption, the claim follows 
from the definition of $U$.
\end{proof}

\subsection{Nilpotent displays}

Let $R$ be a ring which is complete and separated
in the $\Fc$-adic topology for an ideal $\Fc\subset R$
containing $p$. We consider the Witt frame
$$
\WWW_R=(W(R),I_R,R,f,f_1).
$$
Here $I_R\subseteq\Rad R$ as required because
$W(R)=\varprojlim W_n(R/\Fc^n)$ and the successive
kernels in this projective system are nilpotent.
The inverse image of $\Fc$ is a ideal of definition
$J\subset W(R)$.
Nilpotent windows over $\WWW_R$ with respect to $J$
are displays over $R$ which are nilpotent over $R/\Fc$.
By \cite{Zink-Disp} and \cite{Lau-Disp} these are equivalent
to $p$-divisible groups over $R$ which are infinitesimal
over $R/\Fc$.
(Here one uses that displays and $p$-divisible
groups over $R$ are equivalent to compatible systems
of the same objects over $R/\Fc^n$ for $n\geq 1$; cf.\
Lemma \ref{Le-lim} above and \cite[Lemma 4.16]{Messing-Crys}.)

Assume that $R'=R/\Fb$ for a closed ideal $\Fb\subseteq\Fc$ 
equipped with (not necessarily nilpotent) divided powers. 
One can define a factorisation
$$
\WWW_R\xrightarrow{\alpha_1}
\WWW_{R/R'}=(W(R),\tilde I,R',f,\tilde f_1)
\xrightarrow{\alpha_2}\WWW_{R'}
$$
of the projection of frames $\WWW_R\to\WWW_{R'}$ as follows.
Necessarily we define $\tilde I=I_R+W(\Fb)$. The divided
Witt polynomials define an isomorphism
$$
\log:W(\Fb)\cong\Fb^\infty,
$$
and $\tilde f_1:\tilde I_R\to W(R)$ extends 
$f_1$ such that $\tilde f_1[b_0,b_1,\ldots]=[b_1,b_2,\ldots]$
in logarithmic coordinates on $W(\Fb)$. The assumption
$\Fb\subseteq\Fc$ implies that $J$ is an ideal
of definition for $\WWW_{R/R'}$ as well. 

We assume that
the $\Fc$-adic topology of $R$ can be defined by
a sequence of ideals $R\supset I_1\supset I_2\ldots$
such that each $\Fb\cap I_n$ is stable under the
divided powers of $\Fb$. This is automatic when
$\Fc$ is nilpotent or when $R$ is noetherian; 
cf.\ the proof of Proposition \ref{Pr-crys-Dieu}.

\begin{Prop}
\label{Pr-nil-crys-Witt}
The homomorphism $\alpha_2$ is nil-crystalline
with respect to the ideal of definition
$J\subset W(R)$ for both frames.
\end{Prop}

This is essentially \cite[Theorem 44]{Zink-Disp}.

\begin{proof}
By a limit argument the assertion is reduced to the
case where $\Fc\subset R$ is a nilpotent ideal;
see Lemma \ref{Le-lim}. Then Theorem \ref{Th-nil-crys}
applies: The required filtration of $\Fa=W(\Fb)$ is
$\Fa_i=p^i\Fa$. The condition $J^n\Fa=0$ for large $n$
is satisfied because $J^n\subseteq I_R$ for some $n$
and $I_R^{n+1}\subseteq p^nW(R)$ for all $n$,
and $W(\Fb)\cong\Fb^\infty$ is annihilated by some 
power of $p$.
\end{proof}

%----------------------------------------------------------------
\subsection{The main frame}

Let $R$ be a complete regular local ring with 
arbitrary residue field $k$ of characteristic $p$.
Let $C$ be a $p$-ring with residue field $k$.
We choose a surjective ring homomorphism
$$
\FS=C[[x_1,\ldots,x_r]]\to R
$$
that lifts the identity of $k$ such that $x_1,\ldots,x_r$
map to a regular system of parameters for $R$.
There is a power series $E\in\FS$ with constant term $p$
such that $R=\FS/E\FS$. Let $\sigma:C\to C$ be a ring endomorphism
which induces the Frobenius on $\FS/p\FS$ and preserves
the ideal $(x_1,\ldots,x_r)$. We consider the frame
$$
\BBB=(\FS,I,R,\sigma,\sigma_1)
$$
where $\sigma_1(Ey)=\sigma(y)$. Here $\theta=\sigma(E)$.
The proof of Lemma \ref{Le-B-kappa} shows
that $\BBB$ is again a $\varkappa$-frame,
so we have a $u$-homomorphism of frames
$$
\varkappa:\BBB\to\WWW_R.
$$
Let $\Fm\subset\FS$ and $\Fn\subset W(R)$ be the maximal ideals.

\begin{Thm}
\label{Th-kappa-nil-crys}
The homomorphism $\varkappa$ is nil-crystalline
with respect to the ideals of definition 
$\Fm$ of $\BBB$ and $\Fn$ of $\WWW_R$. 
\end{Thm}

\begin{proof}
The proof of Theorem \ref{Th-kappa-crys-gen} applies with 
the following modification: The initial case $a=0$ is
not trivial because $C\not\cong W(k)$ if $k$ is not perfect, 
but one can apply \cite[Theorem 1.6]{Zink-Windows}.
In the diagram \eqref{Diag}, the homomorphisms
$\pi'$ and $\pi$ are only nil-crystalline in general;
whether $\pi$ is crystalline depends on the choice of $\sigma$.
\end{proof}

\subsection{Connected group schemes}

One defines Breuil windows relative to $\FS\to R$ and 
Breuil modules relative to $\FS\to R$ as before. 
A Breuil window $(Q,\phi)$ or a Breuil module $(M,\varphi,\psi)$ 
is called nilpotent if $\phi$ or $\varphi$ is nilpotent
modulo the maximal ideal of $\FS$.
The proof of Lemma \ref{Le-Br-win} shows that
nilpotent Breuil windows are equivalent to
nilpotent $\BBB$-windows. 
Hence Theorem \ref{Th-kappa-nil-crys} implies:

\begin{Cor}
Connected $p$-divisible groups over $R$ are equivalent 
to nilpotent Breuil windows relative to $\FS\to R$.
\qed
\end{Cor}

Similarly we have:

\begin{Thm}
Connected finite flat group schemes over $R$ of
$p$-power order are equivalent to 
nilpotent Breuil modules relative to $\FS\to R$.
\end{Thm}

This is proved like Theorem \ref{Th-modules},
using two additional remarks:

\begin{Lemma}
Every connected finite flat group scheme $H$ over $R$
is the kernel of an isogeny of connected $p$-divisible
groups.
\end{Lemma}

\begin{proof}
We know that $H$ is the kernel of an isogeny 
of $p$-divisible groups $G\to G'$.
There is a functorial exact sequence of $p$-divisible 
groups $0\to G_0\to G\to G_1\to 0$ 
where $G_0$ is connected and $G_1$ 
is etale. Since $\Hom(H,G_1)$ is zero, 
$H$ is the kernel of $G_0\to G_0'$.
\end{proof}

\begin{Lemma}
Every nilpotent Breuil module $(M,\varphi,\psi)$ is the
cokernel of an isogeny of nilpotent Breuil windows.
\end{Lemma}

\begin{proof}
We know that $(M,\varphi,\psi)$ is the cokernel of an 
isogeny of Breuil windows $(Q,\phi)\to(Q',\phi')$. 
There is a functorial exact sequence of 
Breuil windows $0\to Q_0\to Q\to Q_1\to 0$ where $Q_0$
is nilpotent and where $Q_1$ is etale in the sense
that $\phi:Q_1\to Q_1^{(\sigma)}$ is bijective.
Indeed, by \cite[Lemma 10]{Zink-DDisp} it suffices to
construct the sequence over $k$, and $Q_0\otimes_\FS k$ 
is the kernel of $(\phi_k)^n$ for large $n$, where
$\phi_k:Q\otimes_\FS k\to Q^{(\sigma)}\otimes_\FS k$
is the special fibre of $\phi$.

We claim that $Q_1$ and $Q_1'$ have the same rank.
We identify $C$ with $\FS/(x_1,\ldots x_r)$. 
Since $Q\to Q'$ becomes bijective over $\FS[1/p]$, the map
$Q\otimes_\FS C\to Q'\otimes_\FS C$ becomes bijective over $C[1/p]$.
Hence the etale parts $(Q\otimes_\FS C)_1$ and 
$(Q'\otimes_\FS C)_1$ 
have the same rank. This proves the claim because 
$(Q\otimes_\FS C)_1=Q_1\otimes_\FS C$ and similarly for $Q'$.

Let us consider $\bar M=Q_1'/Q_1$. Here $\phi'$ induces 
a homomorphism $\bar\varphi:\bar M\to\bar M^{(\sigma)}$, 
which is surjective as $Q_1'$ is etale.
The natural surjection $\pi:M\to\bar M$  
satisfies $\pi^{(\sigma)}\varphi=\bar\varphi\pi$.
As $\varphi_k$ is nilpotent it follows that $\bar\varphi_k$ 
is nilpotent, thus $\bar M=0$ by Nakayama's lemma.
Hence $Q_1\to Q_1'$ is bijective
because both sides are free of the same rank,
and consequently $M=Q_0'/Q_0$ as desired.
\end{proof}

%----------------------------------------------------------------
%\cleardoublepage
\section{Generalised frames}
\label{Se-gen}

We mention a generalisation of the notion of frames 
and windows which is not considered in the main text.

\begin{Defn}
A generalised frame is a sextuple 
$$\FFF=(S,I,R,\sigma,\sigma_1,\theta)
$$
consisting of a ring $S$, an ideal $I$ of $S$, the quotient ring 
$R=S/I$, a ring endomorphism $\sigma:S\to S$, a $\sigma$-linear homomorphism of $S$-modules $\sigma_1:I\to S$, 
and an element $\theta\in S$, such that we have:
\begin{enumerate}
\renewcommand{\theenumi}{\roman{enumi}}
\item
$I+pS\subseteq\Rad(S)$,
\item
$\sigma(a)\equiv a^p$ mod $pS$ for $a\in S$,
\item
$\sigma(a)=\theta\sigma_1(a)$ for $a\in I$.
\end{enumerate}
\end{Defn}

Since $\sigma_1(I)$ need not generate $S$,
the element $\theta$ need not be determined by
the rest of the data (cf.\ Lemma \ref{Le-theta}).
For a $u$-homomorphism of generalised frames
$\alpha:\FFF\to\FFF'$ we demand that 
$\alpha(\theta)=u\theta'$.

\begin{Defn}
A window $\PPP$ over a generalised frame $\FFF$ is a quadruple
$\PPP=(P,Q,F,F_1)$ where $P$ is a finitely generated
projective $S$-module, $Q\subseteq P$ is a submodule,
$F:P\to P$ and $F_1:Q\to P$ are $\sigma$-linear homomorphisms
of $S$-modules, such that:
\begin{enumerate}
\item 
There is a decomposition $P=L\oplus T$ with $Q=L\oplus IT$,
\item 
$F_1(ax)=\sigma_1(a)F(x)$ for $a\in I$ and $x\in P$,
\item 
$F(x)=\theta F_1(x)$ for $x\in Q$,
\item 
$F_1(Q)+F(P)$ generates $P$ as an $S$-module.
\end{enumerate}
\end{Defn}

If $\FFF$ is a frame this is equivalent to Definition \ref{Def-win}.
The results of sections \ref{Se-fr-win}--\ref{Se-abs-def} hold for
generalised frames as well. Details are left to the interested reader.

%----------------------------------------------------------------

\end{document}